\documentclass[reqno,12pt]{amsart}
\usepackage{amsmath, latexsym, amsfonts, amssymb, amsthm, amscd}
\usepackage{graphicx, color}
\usepackage{pstricks}

\numberwithin{equation}{section}

\newtheorem{theorem}{Theorem}[section]
\newtheorem{lemma}[theorem]{Lemma}
\newtheorem{proposition}[theorem]{Proposition}
\newtheorem{corollary}[theorem]{Corollary}
\newtheorem{rem}[theorem]{Remark}

\newtheorem{example}[theorem]{Example}

\renewcommand{\tilde}{\widetilde}          
\DeclareMathSymbol{\leqslant}{\mathalpha}{AMSa}{"36} 
\DeclareMathSymbol{\geqslant}{\mathalpha}{AMSa}{"3E} 
\DeclareMathSymbol{\eset}{\mathalpha}{AMSb}{"3F}     
\renewcommand{\leq}{\;\leqslant\;}                   
\renewcommand{\geq}{\;\geqslant\;}                   

\newcommand{\R}{\mathbb{R}}

\newcommand{\N}{\mathbb{N}}

\title[Gaussian multiplicative chaos revisited]{Gaussian multiplicative chaos revisited}
\author{}

\begin{document}

\maketitle
\begin{center}
{Raoul Robert \\
\footnotesize \noindent
 Institut Fourier, universit{\'e} Grenoble 1, UMR CNRS 5582, \\ 
 100, rue des Math{\'e}matiques, BP 74, 38402 Saint-Martin d'H{\`e}res cedex, France}

{\footnotesize \noindent e-mail: \texttt{Raoul.Robert@ujf-grenoble.fr}}

\bigskip

{Vincent Vargas \\
\footnotesize 
 
 CNRS, UMR 7534, F-75016 Paris, France \\
  Universit{\'e} Paris-Dauphine, Ceremade, F-75016 Paris, France} \\

{\footnotesize \noindent e-mail: \texttt{vargas@ceremade.dauphine.fr}}
\end{center}

\begin{abstract}
In this article, we extend the theory of multiplicative chaos for positive definite functions in $\R^d$ of the form $f(x)=\lambda^2\ln^+\frac{T}{|x|}+g(x)$ where $g$ is a continuous and bounded function. The construction is simpler and more general than the one defined by Kahane in 1985. As main application, we give a rigorous mathematical meaning to the Kolmogorov-Obukhov model of energy dissipation in a turbulent flow.      

\end{abstract}
\vspace{1cm}
\footnotesize


\noindent{\bf Key words or phrases:} Random measures, Gaussian processes, Multifractal processes.

\noindent{\bf MSC 2000 subject classifications: 60G57, 60G15, 60G25, 28A80}

\normalsize

\section{Introduction}
The theory of multiplicative chaos was first defined rigorously by Kahane in 1985 in the article \cite{cf:Kah}. More specifically, 
Kahane built a theory relying on the notion of $\sigma$-positive type kernel: a generalized function $K:\R^d \times \R^d \rightarrow \R_{+}\cup \lbrace \infty \rbrace$ is of $\sigma$-positive type if there exists a sequence $K_{k}:\R^d \times \R^d \rightarrow \R_{+}$ of continuous positive and positive definite kernels such that:
\begin{equation}\label{eq:defsigma}
K(x,y)=\sum_{k \geq 1}K_{k}(x,y).
\end{equation}
If $K$ is a $\sigma$-positive type kernel with decomposition (\ref{eq:defsigma}), one can consider a sequence of gaussian processes $(X_{n})_{n \geq 1}$ of covariance given by $\sum_{k=1}^{n}K_{k}$. It is proven in \cite{cf:Kah} that the sequence of random measures $m_{n}$ given by:  
\begin{equation}\label{eq:deformelle}
\forall A \in \mathcal{B}(\R^d), \quad m_{n}(A)=\int_{A}e^{X_{n}(x)-\frac{1}{2}E[X_{n}(x)^2]}dx
 \end{equation}  
converges almost surely in the space of Radon measures (equiped with the topology of weak convergence) towards a random measure $m$ and that the limit measure $m$ obtained does not depend on the sequence $(K_{k})_{k \geq 1}$ used in the decomposition (\ref{eq:defsigma}) of $K$. Thus, the theory 
enables to give a unique and mathematically rigorous definition to a random measure $m$ in $\R^d$ defined formally by:
\begin{equation}\label{eq:deformelle}
\forall A \in \mathcal{B}(\R^d), \quad m(A)=\int_{A}e^{X(x)-\frac{1}{2}E[X(x)^2]}dx
 \end{equation}  
where $(X(x))_{x \in \R^d}$ is a "gaussian field" whose covariance $K$ is a $\sigma$-positive type kernel. As it will appear later the $\sigma$-positive type condition is not easy to check in practice. From where the need to avoid this hypothesis.

The main application of the theory is to give a meaning to the "limit-lognormal" model introduced by Mandelbrot in \cite{cf:Man}. In the sequel, we define $\ln^+x$ for $x>0$ by the following formula:
\begin{equation*}
\ln^+x=\max(\ln(x),0).
\end{equation*}
The "limit-lognormal" model corresponds to the choice of a stationnary $K$ given by:
\begin{equation}\label{eq:modeleog}
K(s,t)=\lambda^2 \ln^{+}(R/|x-y|)+O(1) 
\end{equation}
where $\lambda^2,R$ are positive parameters and $O(1)$ is a bounded quantity as $|x-y| \rightarrow 0$. This model has many applications 
that we now review in the following subchapters.

\subsection{Multplicative chaos in dimension 1: a model for the volatility of a financial asset}
If $(X(t))_{t \geq 0}$ is the logarithm of the price of a financial asset, the volatility $m$ of the asset on the interval $[0,t]$ is by definition equal to the quadratic variation of $X$:
\begin{equation*}
m[0,t]=\lim_{n \to \infty}\sum_{k=1}^{n}(X(tk/n)-X(t(k-1)/n))^2
\end{equation*}
The volatility $m$ can be viewed as a random measure on $\R$. The choice for $m$ of multiplicative chaos associated to the kernel $K(s,t)=\lambda^2\ln^{+}\frac{T}{|t-s|}$ satisfies many empirical properties measured on financial markets: lognormality of the volatility, long range correlations (see \cite{cf:Cizeau} for a study of the SP500 index and components and \cite{cf:Co} for a general review). Note that $K$ is indeed of $\sigma$-positive type (see example \ref{ex:logplus} below) so $m$ is well defined. In the context of finance, $\lambda^2$ is called the intermittency parameter in analogy with turbulence and $T$ is the correlation length. Volatility modeling and forecasting is an important field of finance since it is related to option pricing and risk forecasting: we refer to \cite{cf:DuRoVa} for the problem of forecasting volatility with this choice of $m$.   

Given the volatility $m$, the most natural way to construct a model for the (log) price $X$ is to set:
\begin{equation}\label{eq:MRW}
X(t)=B_{m[0,t]}                      
\end{equation}
where $(B_{t})_{t \geq 0}$ is a brownian motion independent of $m$. Formula (\ref{eq:MRW}) defines the Multifractal Random Walk (MRW) first introduced in \cite{cf:BaDeMu} (see \cite{cf:BaKoMu} for a recent review of the financial applications of the MRW model).

\subsection{Multiplicative chaos in dimension 3: a model for the energy dissipation in a turbulent fluid}
We refer to \cite{cf:Fr} for an introduction to the statistical theory of 3 dimensional turbulence. Consider a stationnary flow at high Reynolds number; it is believed that at small scales the velocity field of the flow is homogeneous and isotropic in space . By small scales we mean scales much smaller than the integral scale $R$ charactersistic of the time stationnary force driving the flow. In the work \cite{cf:Kol} and \cite{cf:Obu}, Kolmogorov and Obukhov propose to model the mean energy dissipation per unit mass in a ball $B(x,l)$ of center $x$ and radius $l<<R$ by a random variable $\epsilon_{l}$ such that $\ln(\epsilon_{l})$ is normal  with variance $\sigma_{l}^2$ given by:
\begin{equation*}
\sigma_{l}^2=\lambda^2\ln(\frac{R}{l})+A
\end{equation*}      
where $A$ is a constant and $\lambda^2$ is the intermittency parameter. As noted by Mandelbrot (\cite{cf:Man}), the only way to define such a model is to construct a random measure $\epsilon$ by a limit procedure. Then, one can define $\epsilon_{l}$ by the formula:
\begin{equation*}
 \epsilon_{l}=\frac{3<\epsilon>}{4 \pi l^3}\epsilon(B(x,l))
\end{equation*}
where $<\epsilon>$ is the average mean energy disspation per unit mass. 
Formally, one is looking for a random measure $\epsilon$ such that:
\begin{equation}\label{eq:epsilon}
\forall A \in \mathcal{B}(\R^d), \quad \epsilon(A)=\int_{A}e^{X(x)-\frac{1}{2}E[X(x)^2]}dx
 \end{equation}  
where $(X(x))_{x \in \R^d}$ is a "gaussian field" whose covariance $K$ is given by $K(x,y)=\lambda^2\ln^{+}\frac{R}{|x-y|}$. The kernel $\lambda^2\ln^{+}\frac{R}{|x-y|}$ is positive definite considered as a tempered distribution (see (\ref{eq:positive}) for a definition of positive definite distributions  and lemma \ref{lem:positive} for a proof). Therefore, one can give a rigorous meaning to (\ref{eq:epsilon}) by using theorem-definition \ref{th:chaos} below. 

However, it is not clear if 
$\lambda^2\ln^{+}\frac{R}{|x-y|}$ is of $\sigma$-positive type in $\R^3$ and therefore, in \cite{cf:Kah}, Kahane considers the $\sigma$-positive type kernel $K(x,y)=\int_{1/T}^{\infty}\frac{e^{-u|x-y|}}{u}du$ as an approximation of $\lambda^2\ln^{+}\frac{R}{|x-y|}$: indeed, one can show that $\int_{1/T}^{\infty}\frac{e^{-u|x-y|}}{u}du=\ln^{+}\frac{R}{|x-y|}+g(|x-y|)$ where $g$ is a bounded continuous function. Nevertheless, it is important to work with $\lambda^2\ln^{+}\frac{R}{|x-y|}$ since this choice leads to measures which exhibit generalized scale invariance properties (see proposition \ref{prop:invariance}).           

\subsection{Organization of the paper}
In section 2, we remind the definition of positive definite tempered distributions and we state theorem-definition \ref{th:chaos} where we define multiplicative chaos $m$ associated to kernels of the type $\ln^{+}\frac{R}{|x|}+O(1)$. In section 3, we review the main properties of the measure $m$: existence of moments and density with respect to Lebesgue measure, multifractality and generalized scale invariance.
In section 4 and 5, we give respectively the proofs of section 2 and 3.

\section{Definition of multiplicative chaos}

\subsection{Positive definite tempered distributions}
Let $\mathcal{S}(\R^d)$ be the Schwartz space of smooth rapidly decreasing functions and $\mathcal{S}'(\R^d)$ the space of tempered distributions (see \cite{cf:Sch}). A distribution $f$ in $\mathcal{S}'(\R^d)$ is of positive definite if:
\begin{equation}\label{eq:positive}
\forall \varphi \in \mathcal{S}(\R^d), \quad \int_{\R^d} \int_{\R^d} f(x-y) \varphi(x) \overline{\varphi(y)}dx dy \geq 0.
\end{equation}   
On $\mathcal{S}'(\R^d)$, one can define the Fourier transform $\hat{f}$ of a tempered distribution by the formula:
\begin{equation}\label{eq:fourierpositive}
\forall \varphi \in \mathcal{S}(\R^d), \quad \int_{\R^d}\hat{f}(\xi)\varphi(\xi)d\xi= \int_{\R^d}f(x)\hat{\varphi}(x)dx
\end{equation}   
where $\hat{\varphi}(x)=\int_{\R^d}e^{-2i\pi x.\xi}\varphi(\xi)d\xi$ is the Fourier transform of $\varphi$. An extension of Bochner's theorem (Schwartz, \cite{cf:Sch}) states that a tempered distribution $f$ is  positive definite if and only if it's Fourier transform is a tempered positive measure.

By definition, a function $f$ in $\mathcal{S}'(\R^d)$ is of $\sigma$-positive type if the associated kernel $K(x,y)=f(x-y)$ is of $\sigma$-positive type. As mentioned in the introduction, Kahane's theory of multiplicative chaos is defined for $\sigma$-positive type functions $f$. The main problem stems from the fact that definition (\ref{eq:defsigma}) is not practical. Indeed, is there a simple characterization (like the computation of a Fourier transform) of functions whose associated kernel can be decomposed along (\ref{eq:defsigma})? If such a characterization exists, how does one find the kernels $K_{n}$ explicitely?     

Finally, we recall the following simple implication: If $f$ belongs to $\mathcal{S}'(\R^d)$and is of $\sigma$-positive type, $f$ is positive and positive definite. However, the converse statement is not clear.

\subsection{A generalized theory of multiplicative chaos}
In this subsection, we contruct a theory of multiplicative chaos for positive definite functions of type $\lambda^2\ln^+\frac{R}{|x|}+O(1)$ without the assumption of $\sigma$-positivity for the underlying function. The theory is therefore much easier to use.    

We consider in $\R^d$  a positive definite function $f$ such that 
\begin{equation}\label{eq:defpos}
f(x)=\lambda^2\ln^+\frac{R}{|x|}+g(x)
\end{equation}
where $\lambda^2\not = 2d$ and $g(x)$ is a bounded continuous function.  
Let $\theta: \R^d \rightarrow \R$  be some continuous function with the following properties: 
\begin{enumerate}
\item
$\theta$ is positive definite
\item
$|\theta(x)| \leq \frac{1}{1+|x|^{d+\gamma}}$ for some  $\gamma>0$.
\item
 $\int_{\R^d}\theta(x)dx=1$  
\end{enumerate}

Here is the main theorem of the article:

\begin{theorem}\label{th:chaos}{(Definition of multiplicative chaos)}

For all $\epsilon >0$, we consider the centered gaussian field $(X_{\epsilon}(x))_{x \in \R^d}$ defined by the convolution:
\begin{equation*}
E[X_{\epsilon}(x)X_{\epsilon}(y)]=(\theta^{\epsilon} \ast f)(y-x),
\end{equation*}
where $\theta^{\epsilon}=\frac{1}{\epsilon^d}\theta(\frac{.}{\epsilon})$.
Then the associated random measure $m_{\epsilon}(dx)=e^{X_{\epsilon}(x)-\frac{1}{2}E[X_{\epsilon}(x)^2]}dx$ converges in law in the space of Radon mesures (equiped with the topology of weak convergence) as $\epsilon$ goes to $0$ towards a random measure $m$ independent of the choice of the regularizing function $\theta$ with the properties (1), (2), (3). We call the measure $m$ the multiplicative chaos associated to the function $f$.   
\end{theorem}

We review below two possible choices of the underlying function $f$. The first example is a $d$-dimensional generalization of the cone construction considered in \cite{cf:BaMu}. The second example is $\lambda^2\ln^+\frac{R}{|x|}$ for $d=1,2,3$ (the case $d=2,3$ seems to have never been considered in the litterature). Both examples are in fact of $\sigma$-positive type (except perhaps the crucial example of $\lambda^2\ln^+\frac{R}{|x|}$ in dimension $d=3$) and it is easy to show that in these cases theorem-definition \ref{th:chaos} and Kahane's theory lead to the same limit measure $m$.

\begin{example}\label{ex:Bacry}
One can construct a positive definite function $f$ with decomposition (\ref{eq:defpos}) by generalizing to dimension $d$ the cone construction of \cite{cf:BaMu}. This was performed in \cite{cf:Cha}. For all $x$ in $\R^d$, we define the cone $C(x)$ in $\R^{d} \times \R_{+}$:
\begin{equation*}
C(x)= \lbrace (y,t) \in \R^{d} \times \R_{+}; \; |y-x| \leq \frac{t \wedge T}{2}   \rbrace.
\end{equation*}
The function $f$ is given by:
\begin{equation}\label{eq:Bacry}
f(x)=\lambda^2 \int_{C(0) \cap C(x)}\frac{dydt}{t^{d+1}}
\end{equation}
One can show that $f$ has decomposition (\ref{eq:defpos}) (see \cite{cf:Cha}). The function $f$ is of $\sigma$-positive type in the sense of Kahane since one can write $f=\sum_{n \geq 1}f_{n}$ with $f_n$ given by:
\begin{equation*}
f_{n}(x)=\lambda^2 \int_{C(0) \cap C(x); \; \frac{1}{n} \leq t <  \frac{1}{n-1}}\frac{dydt}{t^{d+1}}.  
 \end{equation*}
In dimension $d=1$, we get the simple formula $f(x)=\ln^+\frac{R}{|x|}$.
\end{example}

\begin{example}\label{ex:logplus}
In dimension $d=1,2$, the function $f(x)=\ln^+\frac{R}{|x|}$ is of $\sigma$-positive type in the sense of Kahane and in particular positive definite. Indeed, one has by straightforward calculations:
\begin{equation*}
\ln^+\frac{T}{|x|}=\int_{0}^{\infty}(t-|x|)_{+}\nu_{T}(dt)
\end{equation*}
where $\nu_{T}(dt)=1_{[0,T[}(t)\frac{dt}{t^2}+\frac{\delta_{T}}{T}$.
For all $\mu>0$, we have:
\begin{equation*}
\ln^+\frac{T}{|x|}=\frac{1}{\mu}\ln^+\frac{T^\mu}{|x|^\mu}=\frac{1}{\mu}\int_{0}^{\infty}(t-|x|^\mu)_{+}\nu_{T^\mu}(dt).
\end{equation*}
We are therefore led to considering the $\mu>0$ such that $(1-|x|^\mu)_{+}$ is positive definite (the so called Kuttner-Golubov problem: see \cite{cf:Gn} for an introduction).

For $d=1$, it is straightforward to show that $(1-|x|)_{+}$ is of positive type. One can thus write $f=\sum_{n \geq 1}f_{n}$ with $f_n$ given by:
\begin{equation*}
f_{n}(x)=\int_{\frac{T}{n}}^{\frac{T}{n-1}}(t-|x|)_{+}\nu_{T}(dt).
 \end{equation*}

For $d=2$, the function $(1-|x|^{1/2})$ is positive definite (Pasenchenko, \cite{cf:PaYu}). One can thus write $f=\sum_{n \geq 1}f_{n}$ with $f_n$ given by:
\begin{equation*}
f_{n}(x)=\int_{\frac{T^{1/2}}{n}}^{\frac{T^{1/2}}{n-1}}(t-|x|^{1/2})_{+}\nu_{T^{1/2}}(dt).
 \end{equation*} 

In dimension $d=3$, the function $\ln^+\frac{R}{|x|}$ is positive definite (see lemma \ref{lem:positive} below) but it is an open question whether it is of $\sigma$-positive type. 
\end{example}

\section{Main properties of multiplicative chaos}
In the sequel, we will consider the structure functions $\zeta_{p}$ defined for all $p$ in $\R$ by:
\begin{equation}\label{eq:zeta} 
\zeta_{p}=(d+\frac{\lambda^2}{2})p-\frac{\lambda^2 p^2}{2}.
\end{equation}

\subsection{Multiplicative chaos is equal to $0$ for $\lambda^2 > 2d$}

The following proposition can be seen as a phase transition and shows that the logarithmic kernel is crucial in the theory of multiplicative chaos:

\begin{proposition}\label{prop:null}
If $\lambda^2 > 2d$, the limit measure is equal to $0$.

\end{proposition}

\subsection{Generalized scale invariance}
In this subsection and the following, in view of proposition \ref{prop:null}, we will suppose that $\lambda^2 < 2d$:

Let $m$ be a homogeneous random measure on $\R^d$. We note $B(0,R)$ the ball of center $0$ and radius $R$ in $\R^d$. We say $m$ has the generalized scale invariance property with integral scale $R>0$ if for all $c$ in $]0,1]$ the following equality in law holds:
\begin{equation}\label{eq:scaleinv}
(m(c A))_{A \subset B(0,R)}\overset{(Law)}{=}e^{\Omega_{c}}(m(A))_{A \subset B(0,R)}
\end{equation}
where $\Omega_{c}$ is a random variable independent from $m$. If $m$ is different from $0$, then it is immediate to prove that $(\Omega_{e^{-t}})_{t \geq 0}$ is a Levy process. In the context of gaussian multiplicative chaos, the process $(\Omega_{e^{-t}})_{t \geq 0}$ will be Brownian motion with drift.  

In order to get scale invariance with integral scale $R$, one can choose $f=\ln^{+}\frac{R}{|.|}$. This is possible if and only if $\ln^{+}\frac{R}{|.|}$ is positive defnite. This motivates the following lemma:

\begin{lemma}\label{lem:positive}
Let $d \geq 1$ be the dimension of the space and $R> 0$ the integral scale. We consider the function 
$f: \R^d \rightarrow \R_{+}$ defined by:
\begin{equation*}
f(x)=\ln^{+} \frac{R}{|x|} .
\end{equation*}
The function $f$ is positive definite if and only if $d \leq 3$. 
\end{lemma}

The above choice of $f$ leads to measures that have the generalized scale invariance property.

\begin{proposition}\label{prop:invariance}
Let $d$ be less or equal to 3 and $m$ be the gaussian multplicative chaos with kernel $\lambda^2 \ln^+\frac{R}{|x|}$. Then $m$ is scale invariant; for all $c$ in $]0,1]$, the following equality holds:
\begin{equation}\label{eq:invariance}
(m(c A))_{A \subset B(0,R)}\overset{(Law)}{=}e^{\Omega_{c}}(m(A))_{A \subset B(0,R)},
\end{equation}
where $\Omega_{c}$ is a gaussian random variable independent of $m$ with mean 
$-(d+\frac{\lambda^2}{2})\ln(1/c)$ and variance $\lambda^2\ln(1/c)$ 
\end{proposition}
The proof of the proposition is straightforward. 

\begin{rem}
It remains an open question to construct homogeneous measures in dimension greater or equal to $4$ which are scale invariant.  
\end{rem}

\subsection{Existence of moments and multifractality}
We remind that we suppose that $\lambda^2 < 2d$: this ensures the existence of $\epsilon>0$ such that $\zeta_{1+\epsilon}>d$. Therefore, there exists a unique $p_{\ast}>1$ such that $\zeta_{p_{\ast}}=d$. 
The following two propositions establish the existence of positive and negative moments for the limit measure. 


\begin{proposition}\label{prop:mompos}{(Positive Moments)}

Let $p$ belong to $]0,p_{\ast}[$ and $m$ be the gaussian multiplicative chaos associated to the function $f$ given by (\ref{eq:defpos}). 
For all bounded $A$ in $\mathcal{B}(\R^d)$, 
\begin{equation}\label{eq:mom}
E[m(A)^p]< \infty
\end{equation}
Let $\theta$ be some function satisfying the conditions (1), (2), (3) of section 2.2. With the notations of theorem \ref{th:chaos}, we consider the random measure $m_{\epsilon}$ associated to $\theta$. We have the following convergence for all bounded $A$ in $\mathcal{B}(\R^d)$:
\begin{equation}\label{eq:convmom}
E[m_{\epsilon}(A)^p] \underset{\epsilon \to 0}{\rightarrow}E[m(A)^p].
\end{equation} 
\end{proposition}

\begin{proposition}\label{prop:momneg}{(Negative Moments)}

Let $p$ belong to $]-\infty,0]$ and $m$ be the gaussian multiplicative chaos associated to the function $f$ given by (\ref{eq:defpos}). 
For all $c>0$, 
\begin{equation}\label{eq:mom}
E[m(B(0,c))^p]< \infty
\end{equation}
Let $\theta$ be some function satisfying the conditions (1), (2), (3) of section 2.2. With the notations of theorem \ref{th:chaos}, we consider the random measure $m_{\epsilon}$ associated to $\theta$. We have the following convergence for all $c>0$:
\begin{equation}\label{eq:convmom}
E[m_{\epsilon}(B(0,c))^p] \underset{\epsilon \to 0}{\rightarrow}E[m(B(0,c))^p].
\end{equation} 
\end{proposition}

The following proposition states the existence of the structure functions. 
\begin{proposition}\label{prop:multi}
Let $p$ belong to $]-\infty,p_{\ast}[$. Let $m$ be the gaussian multiplicative chaos associated to the function $f$ given by (\ref{eq:defpos}). There exists some $C_{p}>0$ (independent of $g$ and $R$ in decomposition  (\ref{eq:defpos}): $C_{p}=C_{p}(\lambda^2)$) such that we have the following multifractal behaviour:
\begin{equation}\label{eq:multi}
E[m([0,c]^d)^p]\underset{c \to 0}{\sim}e^{\frac{p(p-1)g(0)}{2}}C_{p}(\frac{c}{R})^{\zeta_{p}}.
\end{equation} 
\end{proposition}

In the next proposition, we will suppose that $d\leq3$ and that $f(x)=\lambda^2 \ln^+\frac{R}{|x|}$. In this case, we can prove the existence of a $C^{\infty}$ density.

\begin{proposition}\label{prop:density}
Let $d$ be less or equal to 3 and $m$ be the gaussian multiplicative chaos with kernel $\lambda^2 \ln^+\frac{R}{|x|}$.
For all $c<R$, the variable $m(B(0,c))$ has a $C^{\infty}$ density with respect to the Lebesgue measure.
\end{proposition}


\section{Proof of theorem \ref{th:chaos}}
\subsection{A few intermediate lemmas}

In order to prove the theorem, we start  by giving some lemmas we will need in the proof. 

\begin{lemma}\label{lem:approx}
Let $\theta$ be some function on $\R^d$ such that there exists $\gamma,C >0$ with $|\theta(x)| \leq \frac{C}{1+|x|^{d+\gamma}}$. Then we have the following convergence:
\begin{equation}\label{eq:sup}
\sup_{|z|>A}|\int_{\R^d} | \theta(v)| \ln |\frac{z}{z-v}|dv|  \underset{A \to \infty}{\rightarrow} 0.
\end{equation}
\end{lemma}

\proof
We have:
\begin{equation*}
\int_{\R^d} |\theta(v)| \ln |\frac{z}{z-v}|dv =  \int_{|v| \leq \sqrt{|z|}} |\theta(v)| \ln |\frac{z}{z-v}|dv
 +  \int_{|v| > \sqrt{|z|}} |\theta(v)| \ln |\frac{z}{z-v}|dv . 
\end{equation*}

\emph{Consider the first term}:
we have $1-\frac{|v|}{|z|} \leq \frac{|z-v|}{|z|} \leq 1+\frac{|v|}{|z|} $ so that for $|v|\leq \sqrt{|z|}$:
\begin{equation*} 
 1-\frac{1}{\sqrt{|z|}} \leq \frac{|z-v|}{|z|} \leq 1+\frac{1}{\sqrt{|z|}},
\end{equation*}
thus we get $|\ln  \frac{|z-v|}{|z|}| \leq \ln(\frac{1}{ 1-\frac{1}{\sqrt{|z|}}}) \leq \frac{1}{\sqrt{|z|}-1}$. We conclude:
\begin{equation*}
 \int_{|v| \leq \sqrt{|z|}} |\theta(v)| \ln |\frac{z}{z-v}|dv \leq \frac{1}{\sqrt{|z|}-1} \int_{\R^d} |\theta(v)| dv
\end{equation*}

\emph{Consider the second term}:
we have:
\begin{equation*}
 \int_{|v| > \sqrt{|z|}} |\theta(v)| \ln |\frac{z}{z-v}|dv \leq \ln|z| \int_{|v| > \sqrt{|z|}} |\theta(v)|dv+ \int_{|v| > \sqrt{|z|}} |\theta(v)| |\ln |z-v||dv 
\end{equation*}
The first term above is obvious; we decompose the second:
\begin{equation*}
 \int_{|v| > \sqrt{|z|}} |\theta(v)| |\ln |z-v||dv=\int_{\sqrt{|z|} < |v| < |z|+1} |\theta(v)| |\ln |z-v||dv+\int_{|v| \geq |z|+1} |\theta(v)| |\ln |z-v||dv
\end{equation*}
For $|v| \geq |z|+1$, we have $1 \leq |z-v| \leq |z||v|$ and thus
\begin{equation*}
0 \leq \ln|z-v| \leq \ln|z|+\ln|v|
\end{equation*}
which enables to handle the corresponding integral. Let us now estimate the remaining term 
$I=\int_{\sqrt{|z|} < |v| < |z|+1} |\theta(v)| |\ln |z-v||dv$. Applying Cauchy Schwarz gives:
\begin{equation*} 
I \leq (\int_{\sqrt{|z|} < |v| < |z|+1} |\theta(v)|^2dv)^{1/2} (\int_{\sqrt{|z|} < |v| < |z|+1} |\ln |z-v||^2dv)^{1/2},
\end{equation*}
from which we get straightforwardly:
\begin{equation*}
I \leq \frac{ C \ln|z|}{|z|^{d/2+\gamma/2-d/4}} \underset{|z| \to \infty}{\rightarrow} 0. 
\end{equation*}

\qed




We will also use the following lemma:

\begin{lemma} \label{lem:max}
 Let $\lambda$ be a positive number such that $\lambda^2 \not = 2$ and $(X_{i})_{1 \leq i \leq n}$ an i.i.d. sequence of centered gaussian variables with variance 
 $\lambda^2\ln(n)$. For all positive $p$ such that $p<\max(\frac{2}{\lambda^2},1)$, there exists $0<x<1$ such that:
 \begin{equation}\label{eq:max}
  E[\sup_{1 \leq i \leq n}e^{pX_{i}-p\frac{\lambda^2}{2}\ln(n)}]=O(n^{xp})
\end{equation}

\end{lemma}

\proof
By Fubini we get:
\begin{align}
& E[\sup_{1 \leq i \leq n}e^{pX_{i}-p\frac{\lambda^2}{2}\ln(n)}] \nonumber \\
& = \int_{0}^{\infty} P(\sup_{1 \leq i \leq n}e^{pX_{i}-p\frac{\lambda^2}{2}\ln(n)} > v)dv  \nonumber \\
& =  \int_{0}^{\infty} P(\sup_{1 \leq i \leq n}X_{i} > \frac{\ln(v)}{p}+\frac{\lambda^2}{2}\ln(n))dv \nonumber \\ 
& =   \int_{-\infty}^{\infty}pe^{pu} P(\sup_{1 \leq i \leq n}X_{i} > u+\frac{\lambda^2}{2}\ln(n))du \nonumber \\
& \leq 1+ \int_{0}^{\infty}pe^{pu} P(\sup_{1 \leq i \leq n}X_{i} > u+\frac{\lambda^2}{2}\ln(n))du, \label{eq:1} 
\end{align}
where we performed the change of variable $u=\frac{\ln(v)}{p}$ in the above identities.
If we define $\bar{F}(u)=P(X_{1} > u)$ then we have:  
\begin{equation*}
P(\sup_{1 \leq i \leq n}X_{i} > u+\frac{\lambda^2}{2}\ln(n))=1-e^{n\ln(1-\bar{F}(u+\frac{\lambda^2}{2}\ln(n)))}.
\end{equation*}

Let $x$ be some positive number such that $0<x<1$. Using (\ref{eq:1}), we get:
\begin{align}
E[\sup_{1 \leq i \leq n}e^{pX_{i}-p\frac{\lambda^2}{2}\ln(n)}]  & \leq  n^{xp} +p \int_{x\ln(n)}^{\infty}e^{pu}(1-e^{n\ln(1-\bar{F}(u+\frac{\lambda^2}{2} \ln(n) )} ) du \nonumber  \\   
& \leq n^{xp} +p n^{xp} \int_{0}^{\infty}
e^{p\tilde{u}}(1-e^{n\ln(1-\bar{F}(\tilde{u}+(\frac{\lambda^2}{2}+x)\ln(n)))})d\tilde{u} \label{eq:2} 
\end{align}

We have:
\begin{align*}
\bar{F}(\tilde{u}+(\frac{\lambda^2}{2}+x)\ln(n))) & = \frac{1}{\sqrt{2\pi}\lambda\sqrt{\ln(n)}}\int_{\tilde{u}+(\frac{\lambda^2}{2}+x)\ln(n)}^{\infty}e^{-\frac{v^2}{2\lambda^2 \ln(n)}}dv \\
& = \frac{n^{-\frac{(\lambda^2/2+x)^2}{2\lambda^2}}}{\sqrt{2\pi}\lambda\sqrt{\ln(n)}}\int_{\tilde{u}}^{\infty}e^{-(\frac{1}{2}+\frac{x}{\lambda^2})\tilde{v}-\frac{\tilde{v}^2}{2\lambda^2\ln(n)}}d\tilde{v},
\end{align*}
where we performed the change of variable $\tilde{v}=v-(\frac{\lambda^2}{2}+x)\ln(n)$.
Thus, we get:
\begin{align}
 & n^{xp} \int_{0}^{\infty}
e^{p\tilde{u}}(1-e^{n\ln(1-\bar{F}(\tilde{u}+(\frac{\lambda^2}{2}+x)\ln(n)))})d\tilde{u} \nonumber \\
& \leq n^{xp+1} \int_{0}^{\infty}
e^{p\tilde{u}}\bar{F}(\tilde{u}+(\frac{\lambda^2}{2}+x)\ln(n))d\tilde{u}  \nonumber  \\
& \leq \frac{n^{xp+1-\frac{(\lambda^2/2+x)^2}{2\lambda^2}}}{\sqrt{2\pi}\lambda\sqrt{\ln(n)}}\int_{0}^{\infty}e^{p\tilde{u}}(\int_{\tilde{u}}^{\infty}e^{-(\frac{1}{2}+\frac{x}{\lambda^2})\tilde{v}-\frac{\tilde{v}^2}{2\lambda^2\ln(n)}}d\tilde{v})d\tilde{u}  \nonumber   \\
& \leq \frac{n^{xp+1-\frac{(\lambda^2/2+x)^2}{2\lambda^2}}}{p\sqrt{2\pi}\lambda\sqrt{\ln(n)}}\int_{0}^{\infty}e^{p\tilde{v}-(\frac{1}{2}+\frac{x}{\lambda^2})\tilde{v}-\frac{\tilde{v}^2}{2\lambda^2\ln(n)}}d\tilde{v}  \nonumber  \\
& \leq \frac{n^{xp+1-\frac{(\lambda^2/2+x)^2}{2\lambda^2}}}{p\sqrt{2\pi}\lambda\sqrt{\ln(n)}}\int_{-\infty}^{\infty}e^{p\tilde{v}-(\frac{1}{2}+\frac{x}{\lambda^2})\tilde{v}-\frac{\tilde{v}^2}{2\lambda^2\ln(n)}}d\tilde{v} \nonumber  \\
& = \frac{n^{xp+\alpha(x, \lambda ,p)}}{p}  \label{eq:3},
\end{align}
with $\alpha(x, \lambda^2,p)=1-\frac{(\lambda^2/2+x)^2}{2\lambda^2}+(p-\frac{1}{2}-\frac{x}{\lambda^2})^2\frac{\lambda^2}{2}$.  We have by combining (\ref{eq:2}) and (\ref{eq:3}):
\begin{equation*}
E[\sup_{1 \leq i \leq n}e^{pX_{i}-p\frac{\lambda^2}{2}\ln(n)}]   \leq  n^{xp} +n^{xp+\alpha(x, \lambda ,p)},
\end{equation*}
We focus on the case $p \in ]\frac{1}{2}+\frac{1}{\lambda^2},\max(\frac{2}{\lambda^2},1)[$ (This implies inequality (\ref{eq:max}) for $p \leq \frac{1}{2}+\frac{1}{\lambda^2}$ by Holders inequality).

\emph{First case}: $\lambda^2<2$.

Note that $\alpha(1, \lambda^2 ,\frac{2}{\lambda^2})=0$ so if $p<\frac{2}{\lambda^2}$ then there exists $0<x<1$ such that 
$\alpha(x, \lambda^2 ,p)<0$.

\emph{Second case}: $\lambda^2>2$.

Note that $\alpha(1, \lambda^2 ,1)=0$ so if $p<1$ then there exists $0<x<1$ such that 
$\alpha(x, \lambda^2 ,p)<0$.

\qed

\subsection{Proof of theorem \ref{th:chaos}}
For sake of simplicity, we give the proof in the case $d=1$, $R=1$ and the function $f(x)=\lambda^2 \ln^+\frac{1}{|x|}$. This is no restriction; indeed, the proof in the general case is an immediate adaptation of the following proof.

\subsubsection{Uniqueness}
Let  $\alpha \in ]0,1/2[$. We consider $\theta$ and $\tilde{\theta}$ two continuous functions satisfying properties (1), (2) and (3). We note:
\begin{equation*}
m(dt)=e^{X(t)-\frac{1}{2}E[X(t)^2]}dt=\lim_{\epsilon \to 0} \; e^{X_{\epsilon}(t)-\frac{1}{2}E[X_{\epsilon}(t)^2]}dt,
\end{equation*}
where $(X_{\epsilon}(t))_{t \in \R}$ is a gaussian process of covariance $q_{\epsilon}(|t-s|)$ with:
\begin{equation*}
q_{\epsilon}(x)=(\theta^{\epsilon} \ast f)(x)=\lambda^2\int_{\R}\theta(v)\ln^{+}(\frac{1}{|x-\epsilon v|})dv
\end{equation*}
We define similarly the measure $\tilde{m}$, $\tilde{X}_{\epsilon}$ and $\tilde{q}_{\epsilon}$ associated to the function $\tilde{\theta}$. Note that we suppose that the random measures $m_{\epsilon}(dt)=e^{X_{\epsilon}(t)-\frac{1}{2}E[X_{\epsilon}(t)^2]}dt$ and $\tilde{m}_{\epsilon}(dt)=e^{\tilde{X}_{\epsilon}(t)-\frac{1}{2}E[X_{\epsilon}(t)^2]}dt$ converge in law in the space of Radon measures: this is no restriction since the equality $E[m_{\epsilon}(A)]=E[\tilde{m}_{\epsilon}(A)]=|A|$ for all bounded $A$ in $\mathcal{B}(\R)$ implies the measures are tight.

We will show that:
\begin{equation*}
E[m[0,1]^\alpha]=E[\tilde{m}[0,1]^\alpha]
\end{equation*}
for $\alpha$ in the interval $]0,1/2[$.
If we define $Z_{\epsilon}(t)(u)=\sqrt{t}\tilde{X}_{\epsilon}(u)+\sqrt{1-t}X_{\epsilon}(u)$ with $X_{\epsilon}(u)$ and 
$\tilde{X}_{\epsilon}(u)$ independent, we get by using a continuous version of lemma \ref{lem:Kah}:
\begin{equation}\label{eq:integral}
E[\tilde{m}_{\epsilon}[0,1]^\alpha]-E[m_{\epsilon}[0,1]^\alpha]=\frac{\alpha(\alpha-1)}{2}\int_{0}^{1}\varphi_{\epsilon}(t)dt,
\end{equation}
with $\varphi_{\epsilon}(t)$ defined by:
\begin{equation*}
\varphi_{\epsilon}(t)=\int_{[0,1]^2}(\tilde{q}_{\epsilon}(|t_{2}-t_{1}|)-q_{\epsilon}(|t_{2}-t_{1}|)E[\mathcal{X}_{\epsilon}(t,t_{1},t_{2})]dt_{1}dt_{2},
\end{equation*}
where $\mathcal{X}_{\epsilon}(t,t_{1},t_{2})$ is given by:
\begin{equation*}
\mathcal{X}_{\epsilon}(t,t_{1},t_{2})=\frac{e^{Z_{\epsilon}(t)(t_1)+Z_{\epsilon}(t)(t_2)-\frac{1}{2}E[Z_{\epsilon}(t)(t_1)^2]-\frac{1}{2}E[Z_{\epsilon}(t)(t_2)^2]}}{(\int_{0}^{1}e^{Z_{\epsilon}(t)(u)-\frac{1}{2}E[Z_{\epsilon}(t)(u)^2]}du)^{2-\alpha}}.
\end{equation*}
We state and prove the following short lemma we will need in the sequel.
\begin{lemma}\label{lem:inter} 
For $A>0$, we denote $C_{A}^{\epsilon}=\underset{|x| \geq A\epsilon}{\sup}|q_{\epsilon}(x)-\tilde{q}_{\epsilon}(x)|$. We have:
\begin{equation*}
\underset{A \to \infty}{\lim}(\underset{\epsilon \to 0}{\overline{\lim}}C_{A}^{\epsilon})=0.
\end{equation*}
\end{lemma}
\proof
Let $|x|  \geq A\epsilon$. If $|x| \geq 1/2$ then $q_{\epsilon}(x)$ and $\tilde{q}_{\epsilon}$ converge uniformly towards $\ln^{+}\frac{1}{|x|}$ thus $q_{\epsilon}(x)-\tilde{q}_{\epsilon}$ converges uniformly to $0$. If $|x|< 1/2$, we write:
\begin{equation*}
q_{\epsilon}(x)=\ln\frac{1}{\epsilon}+Q(x/\epsilon)+R_{\epsilon}(x),
\end{equation*}
where $Q(x)=\int_{\R}\ln\frac{1}{|x-z|}\theta(z)dz$ and $R_{\epsilon}(x)$ converges uniformly to $0$ (for $|x|<1/2$) as $\epsilon \rightarrow 0$. This follows from straightforward calculations. Applying lemma \ref{lem:approx}, we get that $Q(x)=\ln \frac{1}{|x|}+\Sigma(x)$ with $\Sigma(x) \rightarrow 0$ for $|x| \rightarrow \infty$. Thus $Q(x)-\tilde{Q}(x)$ is a continuous function such that for $|x|  \geq A\epsilon$ and $|x| \leq 1/2$  we have:
\begin{equation*}
|q_{\epsilon}(x)-\tilde{q}_{\epsilon}(x)| \leq \underset{|y| \geq A}{\sup}|Q(y)-\tilde{Q}(y)|+\underset{|x| \leq 1/2}{\sup}|R_{\epsilon}(x)-\tilde{R}_{\epsilon}(x)|
\end{equation*}
The result follows.

\qed

One can decompose expression (\ref{eq:integral}) in the following way:
\begin{equation}\label{eq:integral1}
E[\tilde{m}_{\epsilon}[0,1]^\alpha]-E[m_{\epsilon}[0,1]^\alpha]=\frac{\alpha(\alpha-1)}{2}\int_{0}^{1}\varphi_{\epsilon}^{A}(t)dt+\frac{\alpha(\alpha-1)}{2}\int_{0}^{1}\bar{\varphi}_{\epsilon}^{A}(t)dt
\end{equation}
where:
\begin{equation*}
\varphi_{\epsilon}^{A}(t)=\int_{[0,1]^2,|t_{2}-t_{1}| \leq A\epsilon}(\tilde{q}_{\epsilon}(|t_{2}-t_{1}|)-q_{\epsilon}(|t_{2}-t_{1}|)E[\mathcal{X}_{\epsilon}(t,t_{1},t_{2})]dt_{1}dt_{2}
\end{equation*}
and
\begin{equation*}
\bar{\varphi}_{\epsilon}^{A}(t)=\int_{[0,1]^2,|t_{2}-t_{1}| > A\epsilon}(\tilde{q}_{\epsilon}(|t_{2}-t_{1}|)-q_{\epsilon}(|t_{2}-t_{1}|)E[\mathcal{X}_{\epsilon}(t,t_{1},t_{2})]dt_{1}dt_{2}.
\end{equation*}
With the notations of lemma \ref{lem:inter}, we have:
\begin{align*}
| \bar{\varphi}_{\epsilon}^{A}(t) | & \leq \lambda^2 C_{A}^{\epsilon} \int_{[0,1]^2,|t_{2}-t_{1}| > A\epsilon}E[\mathcal{X}_{\epsilon}(t,t_{1},t_{2})]dt_{1}dt_{2}  \\
& \leq \lambda^2 C_{A}^{\epsilon} \int_{[0,1]^2}E[\mathcal{X}_{\epsilon}(t,t_{1},t_{2})]dt_{1}dt_{2} \\
& = \lambda^2 C_{A}^{\epsilon} E[(\int_{0}^{1}e^{Z_{\epsilon}(t)(u)-\frac{1}{2}E[Z_{\epsilon}(t)(u)^2]}du)^{\alpha}] \\ 
& \leq \lambda^2 C_{A}^{\epsilon}. 
\end{align*}
Thus, taking the limit as $\epsilon$ goes to $0$ in (\ref{eq:integral1}) gives:
\begin{equation*}
\underset{\epsilon \to 0}{\overline{\lim}} |E[\tilde{m}_{\epsilon}[0,1]^\alpha]-E[m_{\epsilon}[0,1]^\alpha]| \leq \frac{\alpha(1-\alpha)}{2}\lambda^2\underset{\epsilon \to 0}{\overline{\lim}}C_{A}^{\epsilon}+ \frac{\alpha(1-\alpha)}{2}\underset{\epsilon \to 0}{\overline{\lim}}\int_{0}^{1} | \varphi_{\epsilon}^{A}(t) |dt
\end{equation*}

We will show that $\underset{\epsilon \to 0}{\lim}\varphi_{\epsilon}^{A}(0)=0$ (the general case $\varphi_{\epsilon}^{A}(t)$ is similar). There exists a constant $\tilde{C}_{A}>0$ independent of $\epsilon$ such that:
\begin{equation*}
\underset{|x| \leq A\epsilon}{\sup}|\tilde{q}_{\epsilon}(x)-q_{\epsilon}(x)| \leq \tilde{C}_{A}.
\end{equation*}
Therefore, we have:
\begin{align}
|\varphi_{\epsilon}^{A}(0)| & \leq \tilde{C}_{A} \int_{0}^{1} \int_{t_{1}-A\epsilon}^{t_{1}+A\epsilon}E[\mathcal{X}_{\epsilon}(0,t_{1},t_{2})]dt_{2} dt_{1} \nonumber \\
& = \tilde{C}_{A} E \Big[ \frac{\int_{0}^1\int_{t_1-A\epsilon}^{t_1+A\epsilon}e^{X_{\epsilon}(t_1)+X_{\epsilon}(t_2)-\frac{1}{2}E[X_{\epsilon}(t_1)^2]-\frac{1}{2}E[X_{\epsilon}(t_2)^2]}dt_1dt_2}{(\int_{0}^{1}e^{X_{\epsilon}(u)-\frac{1}{2}E[X_{\epsilon}(u)^2]}du)^{2-\alpha}} \Big]  \label{eq:varphi1} 
\end{align}

Now we have:
\begin{align*}
& \int_{0}^1\int_{t_1-A\epsilon}^{t_1+A\epsilon}e^{X_{\epsilon}(t_1)+X_{\epsilon}(t_2)-\frac{1}{2}E[X_{\epsilon}(t_1)^2]-\frac{1}{2}E[X_{\epsilon}(t_2)^2]}dt_1dt_2 \\
& \leq (\underset{t_1}{\sup}\int_{t_1-A\epsilon}^{t_1+A\epsilon}e^{X_{\epsilon}(t_2)-\frac{1}{2}E[X_{\epsilon}(t_2)^2]}dt_2 ) \int_{0}^1 e^{X_{\epsilon}(t_1)-\frac{1}{2}E[X_{\epsilon}(t_1)^2]}dt_1   \\
& \leq 2  (\underset{0 \leq i < \frac{1}{2A\epsilon}}{\sup}\int_{2iA\epsilon}^{2(i+1)A\epsilon}e^{X_{\epsilon}(t_2)-\frac{1}{2}E[X_{\epsilon}(t_2)^2]}dt_2 ) \int_{0}^1 e^{X_{\epsilon}(t_1)-\frac{1}{2}E[X_{\epsilon}(t_1)^2]}dt_1
\end{align*}
In view of (\ref{eq:varphi1}), this implies:
\begin{align*}
|\varphi_{\epsilon}^{A}(0)| & \leq 2 \tilde{C}_{A} E \Big[  (\underset{0 \leq i < \frac{1}{2A\epsilon}}{\sup}\int_{2iA\epsilon}^{2(i+1)A\epsilon}e^{X_{\epsilon}(t_2)-\frac{1}{2}E[X_{\epsilon}(t_2)^2]}dt_2 ) (\int_{0}^1 e^{X_{\epsilon}(t_1)-\frac{1}{2}E[X_{\epsilon}(t_1)^2]}dt_1)^{\alpha-1}  \Big]  \\
& \leq 2 \tilde{C}_{A} E \Big[  (\underset{0 \leq i < \frac{1}{2A\epsilon}}{\sup}\int_{2iA\epsilon}^{2(i+1)A\epsilon}e^{X_{\epsilon}(t_2)-\frac{1}{2}E[X_{\epsilon}(t_2)^2]}dt_2 )^{\alpha} \Big], 
\end{align*}
where we used the inequality $\frac{\sup_{i}a_{i}}{(\sum_{i}a_i)^{1-\alpha}}\leq (\sup_{i}a_{i})^\alpha$.
For sake of simplicity, we now replace $2A$ by $A$.

The idea to study the above supremum is to make the approximation $X_{\epsilon}(t) \approx X_{\epsilon}(A i \epsilon)$ for $t$ in $[A i \epsilon,A(i+1)\epsilon]$. If we define $\mathcal{C}_{\epsilon}$ by:
\begin{equation}\label{eq:defT} 
\mathcal{C}_{\epsilon}= \underset{\underset{Ai\epsilon \leq u \leq A(i+1)\epsilon}{0 \leq i < \frac{1}{A\epsilon}}}{\sup}(X_{\epsilon}(u)-X_{\epsilon}(Ai \epsilon)),
\end{equation}
then we have:
\begin{align}
& E \Big[ (\sup_{0 \leq i < \frac{1}{A\epsilon}}\int_{Ai\epsilon}^{A(i+1)\epsilon}e^{X_{\epsilon}(t)-\frac{1}{2}E[X_{\epsilon}(t)^2]}dt)^{\alpha} \Big]  \nonumber  \\ 
& \leq  E \Big[ (\sup_{0 \leq i < \frac{1}{A\epsilon}}\int_{Ai\epsilon}^{A(i+1)\epsilon}e^{X_{\epsilon}(A i \epsilon)-\frac{1}{2}E[X_{\epsilon}(A i \epsilon)^2]}dt)^{\alpha} 
e^{\alpha \mathcal{C}_{\epsilon}}\Big]   \nonumber  \\  
& =  E \Big[ (\epsilon A \sup_{0 \leq i < \frac{1}{A\epsilon}} e^{X_{\epsilon}(A i \epsilon)-\frac{1}{2}E[X_{\epsilon}(A i \epsilon)^2]})^{\alpha} 
e^{\alpha \mathcal{C}_{\epsilon}}\Big] \nonumber  \\  
 & \leq (\epsilon A)^{\alpha} E \Big[ ( \sup_{0 \leq i < \frac{1}{A\epsilon}} e^{X_{\epsilon}(A i \epsilon)-\frac{1}{2}E[X_{\epsilon}(A i \epsilon)^2]})^{2\alpha} \Big]^{1/2} 
E \Big[ e^{2\alpha \mathcal{C}_{\epsilon}} \Big]^{1/2}. \label{eq:varphi2} 
 \end{align}

It is straightforward to see that there exists some $c \geq 0$ (independent of $\epsilon$) such that for all $s,t$ in $[0,1]$:
\begin{equation*}
E[X_{\epsilon}(s)X_{\epsilon}(t)] \geq -c
\end{equation*}
We introduce a centered gaussian random variable $Z$ independent of $X_{\epsilon}$ and such that $E[Z^2]=c$.
Let $(R_{i}^{\epsilon})_{1 \leq i < \frac{1}{A \epsilon}}$ be a sequence of i.i.d gaussian random variables such that $E[(R_{i}^{\epsilon})^2]=E[X_{\epsilon}(A i \epsilon)^2]+c$. By applying corollary \ref{cor:Kah}, we get:
\begin{align*}
E \Big[ ( \sup_{0 \leq i < \frac{1}{A\epsilon}} e^{X_{\epsilon}(A i \epsilon)-\frac{1}{2}E[X_{\epsilon}(A i \epsilon)^2]})^{2\alpha} \Big] & = \frac{1}{e^{2\alpha^2c-\alpha c}} E \Big[ ( \sup_{0 \leq i < \frac{1}{A\epsilon}} e^{X_{\epsilon}(A i \epsilon)+Z-\frac{1}{2}E[X_{\epsilon}(A i \epsilon)^2]-\frac{c}{2}})^{2\alpha} \Big]   \\
& \leq  \frac{1}{e^{2\alpha^2c-\alpha c}} E \Big[ ( \sup_{0 \leq i < \frac{1}{A\epsilon}} e^{R_{i}^{\epsilon}-\frac{1}{2}E[(R_{i}^{\epsilon})^2]})^{2\alpha} \Big]  \\
\end{align*} 
We have $E[(R_{i}^{\epsilon})^2]=\lambda^2\ln\frac{1}{\epsilon}+C(\epsilon)$ with $C(\epsilon)$ converging to some constant as $\epsilon$ goes to $0$. Since $2\alpha<1$, by appling lemma \ref{lem:max}, there exists 
$0<x<1$ such that:
\begin{equation*}
 E \Big[ ( \sup_{0 \leq i < \frac{1}{A\epsilon}} e^{R_{i}^{\epsilon}-\frac{1}{2}E[(R_{i}^{\epsilon})^2]})^{2\alpha} \Big] 
 \leq C (\frac{1}{\epsilon})^{2\alpha x}
\end{equation*}
and therefore we have:
\begin{equation*}
|\varphi_{\epsilon}^{A}(0)|  \leq C \epsilon^{\gamma} E \Big[ e^{2\alpha \mathcal{C}_{\epsilon}} \Big]^{1/2}
\end{equation*}
with $\gamma=\alpha(1-x)>0$.


One can write $\mathcal{C}_{\epsilon}=\underset{\underset{0 \leq v \leq 1}{0 \leq i < \frac{1}{A\epsilon}}}{\sup}W_{\epsilon}^{i}(v)$ where $W_{\epsilon}^{i}(v)=X_{\epsilon}(Ai \epsilon+A\epsilon v)-X_{\epsilon}(Ai \epsilon)$.
We have:
\begin{equation*}
E[W_{\epsilon}^{i}(v)W_{\epsilon}^{i}(v')]=g_{\epsilon}(v-v')
\end{equation*}
where $g_{\epsilon}$ is a continuous function bounded by $M$ independently of $\epsilon$. Let $Y$ be a centered Gaussian random variable independent of $W_{\epsilon}^{i}$ such that: $E[Y^2]=M$. Thus, we can write:
\begin{equation*} 
E \Big[ e^{2\alpha \mathcal{C}_{\epsilon}} \Big]=\frac{E \Big[ e^{2\alpha\underset{i,v}{\sup}W_{\epsilon}^{i}(v)}\Big] }{e^{2\alpha^2M^2}}.
\end{equation*}

Now let us consider a family $(\overline{W}_{\epsilon}^{i})_{1 \leq i < \frac{1}{A\epsilon}}$ of centered i.i.d. Gaussian processes of law $(W_{\epsilon}^{0}(v)+Y)_{0 \leq v \leq 1}$. Applying corollary \ref{cor:Kah} of the appendix, we get:
\begin{equation*}
E\Big[ e^{2\alpha  \mathcal{C}_{\epsilon}} \Big]  \leq  \frac{E \Big[ e^{2\alpha\underset{i,v}{\sup}\overline{W}_{\epsilon}^{i}(v)}\Big] }{e^{2\alpha^2M^2}} 
\end{equation*}
We now estimate $E \Big[ e^{2\alpha\underset{i,v}{\sup}\overline{W}_{\epsilon}^{i}(v)}\Big]$. Let us denote $\mathcal{X}_{i}=\underset{0 \leq v \leq 1}{\sup}\overline{W}_{\epsilon}^{i}$. Applying Corollary 3.2 of \cite{cf:LeTa} to the continuous gaussian process $(W_{\epsilon}^{0}(v)+Y)_{0 \leq v \leq 1}$, we get that the random variable has a Gaussian tail:
\begin{equation*}
P(\mathcal{X}_{i}>z) \leq Ce^{-\frac{z^2}{2\sigma^2}}, \quad \forall z >0
\end{equation*}
for some $C$ and $\sigma$. The above tail inequality gives the existence of some constant $C>0$ such that:
\begin{equation*}
E\Big[ e^{2\alpha \underset{0 \leq i < \frac{1}{A\epsilon}}{\sup}\mathcal{X}_{i}} \Big]  \leq Ce^{C\sqrt{\ln(\frac{1} {\epsilon})}}.
\end{equation*}
Therefore we have $E \Big[ e^{2\alpha \mathcal{C}_{\epsilon}} \Big] \leq Ce^{C\sqrt{\ln(\frac{1} {\epsilon})}}$ and then:
\begin{equation*}
 |\varphi_{\epsilon}^{A}(0)| \leq C \epsilon^{\gamma}e^{C\sqrt{\ln(\frac{1} {\epsilon})}}.
\end{equation*}
It follows that $\underset{\epsilon \to 0}{\overline{\lim}} |\varphi_{\epsilon}^{A}(0)|=0$ so that for $\alpha < 1/2$:
\begin{equation*}
\underset{\epsilon \to 0}{\overline{\lim}} |E[\tilde{m}_{\epsilon}[0,1]^\alpha]-E[m_{\epsilon}[0,1]^\alpha]| \leq \frac{\alpha(1-\alpha)}{2}\lambda^2\underset{\epsilon \to 0}{\overline{\lim}}C_{A}^{\epsilon}.
\end{equation*}
Since $\underset{\epsilon \to 0}{\overline{\lim}}C_{A}^{\epsilon}\rightarrow 0$ as $A$ goes to infinity (lemma \ref{lem:inter}), we conclude that:  
\begin{equation*}
\underset{\epsilon \to 0}{\overline{\lim}} |E[\tilde{m}_{\epsilon}[0,1]^\alpha]-E[m_{\epsilon}[0,1]^\alpha]|=0.
\end{equation*}
It is straightforward to check that the above proof can be generalized to show that for all positive $\lambda_{1}, \ldots, \lambda_{n}$ and intervals $I_{1}, \ldots, I_{n}$ we have:
\begin{equation*}
E[(\sum_{k=1}^n\lambda_{k}m(I_{k}))^\alpha]=E[(\sum_{k=1}^n\lambda_{k}\tilde{m}(I_{k}))^\alpha]
\end{equation*}
This implies that the random measures $m$ and $\tilde{m}$ are equal (see \cite{cf:DaVe}).

\subsubsection{Existence}
Let $f(x)$ be a real positive definite function on $\R^d$ (note that this implies that $f$ is symmetric). Let us recall that the centered Gaussian field of correlation $f(x-y)$ is given by:
\begin{equation*}
 X(x)=\int_{\R^d}\zeta(x,\xi)\sqrt{\hat{f}(\xi)}W(d\xi),
\end{equation*}
where $\zeta(x,\xi)=\text{cos}(2\pi x.\xi)-\text{sin}(2\pi x.\xi)$ and $W(d\xi)$ is the standard white noise on $\R^d$. This can also be written:
\begin{equation}\label{eq:cons}
 X(x)=\int_{]0,\infty[ \times \R^d}\zeta(x,\xi)\sqrt{\hat{f}(\xi)}g(t,\xi)W(dt,d\xi),
\end{equation}
where $W(dt,d\xi)$ is the white noise on $]0,\infty[ \times \R^d$ and $\int_{0}^{\infty}g(t,\xi)^2dt=1$ for all $\xi$. The interest of the expression (\ref{eq:cons}) appears in what follows. Let the function $\theta$ be radially symmetric and $\hat{\theta}$ be a decreasing function of $|\xi|$ (for instance take $\theta(x)=\frac{e^{-|x|^2/2}}{(2 \pi)^{d/2}}$). Let us consider $g(t,\xi)=\sqrt{-\hat{\theta}'(t|\xi|)|\xi|}$ so that $\int_{\epsilon}^{\infty}g(t,\xi)^2dt=\hat{\theta}(\epsilon|\xi|)$ for $|\xi| \not = 0$. Then if we consider the fields $X_{\epsilon}$ defined by:
\begin{equation}\label{eq:mart}
X_{\epsilon}(x)=\int_{]\epsilon,\infty[ \times \R^d}\zeta(x,\xi)\sqrt{\hat{f}(\xi)}g(t,\xi)W(dt,d\xi)
\end{equation}
we will find:
\begin{align*}
E[X_{\epsilon}(x)X_{\epsilon}(y)] & = \int_{\R^d} \text{cos}(2\pi(x-y).\xi)\hat{f}(\xi)\hat{\theta}(\epsilon|\xi|)d\xi  \\
& = (f \ast \theta^\epsilon)(x-y).
\end{align*}
The interest of (\ref{eq:mart}) is to make the approximation process appear as a martingale. Indeed, if we define the filtration $\mathcal{F}_{\epsilon}=\sigma \lbrace W(A,x), A \subset ]\epsilon,\infty[, x \in \R^d  \rbrace$, we have that for all $A \in \mathcal{B}(\R^d)$, $(m_{\epsilon}(A))_{\epsilon > 0}$ is a positive $\mathcal{F}_{\epsilon}$-martingale of expectation $|A|$ so it converges almost surely to a random variable $m(A)$ such that:
\begin{equation}\label{eq:major}
E[m(A)] \leq |A|. 
\end{equation}
This defines a collection $(m(A))_{A \in \mathcal{B}(\R^d)}$ of random variables such that:
\begin{enumerate}
\item
for all disjoint and bounded sets $A_{1}, A_{2}$ in $\mathcal{B}(\R^d)$,
\begin{equation*} 
m(A_{1} \cup A_{2})=m(A_{1})+m(A_{2}) \quad a.s.
\end{equation*}
\item
For any bounded sequence $(A_{n})_{n \geq 1}$ decreasing to $\emptyset$:
\begin{equation*} 
m(A_{n})\underset{n \to \infty}{\longrightarrow} 0 \quad a.s.
\end{equation*}

 \end{enumerate}
By theorem 6.1. VI. in \cite{cf:DaVe}, one can consider a version of the collection $(m(A))_{A \in \mathcal{B}(\R^d)}$ such that $m$ is a random measure. It is straightforward that $m_{\epsilon}$ converges almost surely towards $m$ in the space of Radon measures (equiped with the weak topology).

\section{Proofs of section 3}

\subsection{Proof of proposition \ref{prop:null}}

\proof
Since $\zeta_{1}=d$, note that $\lambda^2 > 2d$ is equivalent to the existence of $\alpha<1$ such that $\zeta_{\alpha}>d$. Let $\alpha$ be fixed and such that $\zeta_{\alpha}>d$. We will show that $m[[0,1]^d]=0$. We partition the cube $[0,1]^d$ into $\frac{1}{\epsilon^d}$ subcubes $(I_{j})_{1 \leq j \leq \frac{1}{\epsilon^d}}$ of size $\epsilon$. One has by subadditivity and homogeneity:
\begin{align*}
&E[(\int_{[0,1]^d}e^{X_{\epsilon}(x)-\frac{1}{2}E[X_{\epsilon}(x)^2]}dx)^\alpha] \\
 &=  E[(\sum_{1 \leq j \leq \frac{1}{\epsilon^d}}\int_{I_{j}}e^{X_{\epsilon}(x)-\frac{1}{2}E[X_{\epsilon}(x)^2]}dx)^\alpha]  \\
& \leq  E[\sum_{1 \leq j \leq \frac{1}{\epsilon^d}}(\int_{I_{j}}e^{X_{\epsilon}(x)-\frac{1}{2}E[X_{\epsilon}(x)^2]}dx)^\alpha]  \\
&=\frac{1}{\epsilon^d}E[(\int_{[0,\epsilon]^d}e^{X_{\epsilon}(x)-\frac{1}{2}E[X_{\epsilon}(x)^2]}dx)^\alpha] 
\end{align*}
 
 Let $Y_{\epsilon}$ be a centered gaussian random variable of variance $\lambda^2\ln(\frac{1}{\epsilon})+\lambda^2 c$ where $c$ is such that:
\begin{equation*}
\theta^\epsilon \ast \ln^+\frac{1}{|x|} \geq \ln\frac{1}{\epsilon}+c
\end{equation*}
for $|x| \leq \epsilon$ and $\epsilon$ small enough. By definition of $c$, we have 
\begin{equation*}
\forall t,t' \in [0,\epsilon],  \quad E[X_{\epsilon}(t)X_{\epsilon}(t')] \geq E[Y_{\epsilon}^2].
\end{equation*}
Using corollary (\ref{cor:Kahane}) in a continuous version, this implies:
\begin{align*}
&E[(\int_{[0,1]^d}e^{X_{\epsilon}(t)-\frac{1}{2}E[X_{\epsilon}(t)^2]}dt)^\alpha] \\
& \leq \frac{1}{\epsilon^d}E[(\int_{[0,\epsilon]^d}e^{Y_{\epsilon}-\frac{1}{2}E[Y_{\epsilon}^2]}dt)^\alpha] \\
& = e^{\frac{\alpha^2-\alpha}{2}c}\epsilon^{\zeta_{\alpha}-d}
\end{align*}
Taking the limit as $\epsilon$ goes to $0$ gives $m[[0,1]^d]=0$.
\qed

\subsection{Proof of lemma \ref{lem:positive}}

\proof

One has the following general formula for the Fourier transform of radial functions:
\begin{equation}\label{eq:Schwartz}
\hat{f}(\xi)=\frac{2\pi}{|\xi|^{\frac{d-2}{2}}}\int_{0}^{\infty}\rho^{\frac{d}{2}}J_{\frac{d-2}{2}}(2\pi |\xi| \rho)f(\rho) d\rho,
\end{equation}
where $J_{\nu}$ is the Bessel function of order $\nu$.

\emph{First case}: $d \leq 3$. 

It is enough to consider the case $d=3$ (Indeed, this implies that the same holds for smaller dimensions). Using the explicit formula $J_{\frac{1}{2}}(x)=\sqrt{\frac{2}{\pi x}}\sin(x)$, we conclude by integrating by parts:
\begin{align*}
\hat{f}(\xi) & =\frac{2}{|\xi|}\int_{0}^{T}\rho \sin(2\pi |\xi| \rho)\ln(\frac{T}{\rho}) d\rho  \\
& = \frac{1}{\pi |\xi|^{2}} \int_{0}^{T}\cos(2\pi |\xi| \rho)(\ln(\frac{T}{\rho})-1) d\rho  \\ 
& =  \frac{1}{2 \pi^2 |\xi|^{3}}( \int_{0}^{T} \frac{\sin(2\pi |\xi| \rho)}{\rho}d\rho-\sin(2\pi |\xi| T))  \\
& = \frac{1}{2 \pi^2 |\xi|^{3}}( \text{sinc}(2\pi |\xi| T)-\sin(2\pi |\xi| T)),
\end{align*}
where $ \text{sinc}$ is the sinus cardinal function:
\begin{equation*}
\text{sinc}(x)=\int_{0}^x \frac{\sin(\rho)}{\rho} d \rho.
\end{equation*}
We introduce for $x \geq 0$ the function $l(x)=\text{sinc}(x)-\sin(x)$. We have $l'(x)=\frac{\sin(x)-x\cos(x)}{x}$. Thus, there exists $\alpha$ in $]\pi,2\pi[$ such that  $l$ is increasing on $]0,\alpha[$ and decreasing on $]\alpha,2\pi[$. Since $l(0)=0$ and $l(2\pi)=\int_{0}^{2\pi} \frac{\sin(\rho)}{\rho} d \rho \geq 0$, we conclude that for all $x$ in $[0,2\pi]$, $l(x)\geq 0$. A classical computation (Dirichlet integral) gives $\int_{0}^{\infty} \frac{\sin(\rho)}{\rho} d \rho=\frac{\pi}{2}$. Thus, we have by an integration by parts:
\begin{align*}
\int_{0}^{2\pi} \frac{\sin(\rho)}{\rho} d \rho & =  \frac{\pi}{2}-\int_{2\pi}^{\infty} \frac{\sin(\rho)}{\rho} d \rho  \\
& =  \frac{\pi}{2}-\int_{2\pi}^{\infty} \frac{1-\cos(\rho)}{\rho^2} d \rho  \\
& \geq  \frac{\pi}{2}-\frac{1}{2\pi} \\
& \geq 1
\end{align*}
Therefore, if $x \geq 2\pi$, we have:
\begin{align*} 
l(x) & = \int_{0}^{x} \frac{\sin(\rho)}{\rho} d \rho -\sin(x)  \\
& \geq \int_{0}^{2\pi} \frac{\sin(\rho)}{\rho} d \rho -\sin(x) \\
& \geq 0.
\end{align*}

\emph{Second case}: $d \geq 4$.
Combining (\ref{eq:Schwartz}) with the identity $\frac{d}{dx}(x^{\nu}J_{\nu}(x))=x^{\nu}J_{\nu-1}(x)$, we get:
\begin{align}
\hat{f}(\xi) & =\frac{2\pi}{|\xi|^{\frac{d-2}{2}}}\int_{0}^{T}\rho^{\frac{d}{2}}J_{\frac{d-2}{2}}(2\pi |\xi| \rho)\ln(\frac{T}{\rho}) d\rho  \nonumber \\
& = \frac{1}{(2\pi)^{d/2}|\xi|^{d}}\int_{0}^{2\pi |\xi| T}x^{\frac{d}{2}}J_{\frac{d-2}{2}}(x)\ln(\frac{2\pi |\xi| T}{x}) dx  \nonumber  \\
& =   \frac{1}{(2\pi)^{d/2}|\xi|^{d}}\int_{0}^{2\pi |\xi| T}x^{\frac{d}{2}-1}J_{\frac{d}{2}}(x)dx \label{eq:Bessel1}
\end{align}

One has the following asymptotic expansion as $x$ goes to $\infty$ (\cite{cf:GrMacMa}):
\begin{equation}\label{eq:Bessel2}
 J_{\nu}(x)=\sqrt{\frac{2}{\pi x}}\cos(x-\frac{(1+2\nu)\pi}{4})-\frac{(4\nu^2-1)\sqrt{2}}{8\sqrt{\pi}x^{3/2}}\sin(x-\frac{(1+2\nu)\pi}{4})+O(\frac{1}{x^{5/2}}).
\end{equation}
Combining (\ref{eq:Bessel1}) with (\ref{eq:Bessel2}), we therefore get the following expansion as $|\xi|$ goes to infinity:
\begin{equation*}
\hat{f}(\xi)=\frac{1}{(2\pi)^{d/2}|\xi|^{d}}\left(\sqrt{\frac{2}{\pi}}(2\pi|\xi|T)^{\frac{d-3}{2}}\sin(2\pi|\xi|T-\frac{(1+2\nu)\pi}{4})+o(|\xi|^{\frac{d-3}{2}})\right).
\end{equation*}  
Thus $\underset{|\xi| \to \infty}{\overline{\lim}}|\xi|^{d}\hat{f}(\xi)=-\underset{|\xi| \to \infty}{\underline{\lim}}|\xi|^{d}\hat{f}(\xi)=+\infty$.
\qed 

\subsection{Proofs of section 3.3}

\subsubsection{Proof of proposition \ref{prop:mompos} and proposition \ref{prop:momneg}}
Let $\theta$ be some function satisfying the conditions (1), (2), (3) of section 2.2 and $m_{\epsilon}$ be the random measure associated to $\theta^{\epsilon} \ast f$. We consider $\tilde{m}_{\epsilon}$ the random measure associated to $\tilde{f}_{\epsilon}$ where $\tilde{f}_{\epsilon}$ is the function of example \ref{ex:Bacry}:
\begin{equation*}
\tilde{f}_{\epsilon}(x)=\lambda^2 \int_{C(0) \cap C(x); \; \epsilon < t <  \infty}\frac{dydt}{t^{d+1}}.  
 \end{equation*}
 One can show that there exists $c,C>0$ such that for all $x$ we have:
 \begin{equation*}
\tilde{f}_{\epsilon}(x)-c \leq (\theta^{\epsilon} \ast f)(x) \leq  \tilde{f}_{\epsilon}(x)+C
 \end{equation*}
 By using corollary \ref{cor:Kahane} of the appendix in a continuous version, we conclude that there exists $c,C>0$ such that for all $\epsilon$ and all bounded $A$ in $\mathcal{B}(\R^d)$:
 \begin{equation*}
 c E[\tilde{m}_{\epsilon}(A)^p] \leq E[m_{\epsilon}(A)^p] \leq CE[\tilde{m}_{\epsilon}(A)^p]. 
\end{equation*}
 
 \emph{First case:} $p$ belongs to $]0,p_{\ast}[$.
 
 Proposition  \ref{prop:mompos} is therefore established if we can show that: 
 \begin{equation*}
 \sup_{\epsilon>0}E[\tilde{m}_{\epsilon}(A)^p] < \infty.
\end{equation*}
 The above bound can be proved by adapting the proof of theorem 3 in \cite{cf:BaMu} .
 
 \emph{Second case:} $p$ belongs to $]-\infty,0[$. 
 
 Proposition \ref{prop:mompos} is therefore established if we can show that for all $c>0$: 
 \begin{equation*}
 \sup_{\epsilon>0}E[\tilde{m}_{\epsilon}(B(0,c))^p] < \infty.
\end{equation*}
 The above bound can be proved by adapting the proof of the corresponding result in \cite{cf:Mol}.
 
\subsubsection{Proof of proposition \ref{prop:multi}} 
For the sake of simplicity, we consider the case $R=1$ and we will consider the case $p \in [1,p_{\ast}[$.
 We consider $\theta$ a continuous and positive function with compact support $B(0,A)$ satisfying properties (1), (2) and (3) of section 2.2. We note:
\begin{equation*}
m_{\epsilon}(dx)=e^{X_{\epsilon}(x)-\frac{1}{2}E[X_{\epsilon}(x)^2]}dx,
\end{equation*}
where $(X_{\epsilon}(x))_{x \in \R^d}$ is a gaussian field of covariance $q_{\epsilon}(x-y)$ with:
\begin{equation*}
q_{\epsilon}(x)=(\theta^{\epsilon} \ast f)(x)=\int_{\R^d}\theta(z)(\lambda^2\ln^{+}\frac{1}{|x-\epsilon z|}+g(x-\epsilon z))dz.
\end{equation*}
Let $c,c'$ be two positive numbers in $]0,\frac{1}{2}[$ such that $c<c'$. If $\epsilon$ is sufficiently small and $u,v$ belong to $[0,1]^d$, we get:
\begin{align*}
q_{c\epsilon}(c(v-u)) & = \int_{\R^d}\theta(z) \left( \lambda^2 \ln \frac{1}{|c(v-u)-c \epsilon z|}+g(c(v-u)- c \epsilon z) \right) dz  \\
& =  \lambda^2 \ln(\frac{c'}{c})+ \int_{\R^d}\theta(z) \left( \lambda^2 \ln \frac{1}{|c'(v-u)-c' \epsilon z|}+g(c(v-u)- c \epsilon z) \right)dz  \\
& \leq  \lambda^2 \ln(\frac{c'}{c})+ q_{c'\epsilon}(c'(v-u))+C_{c,c',\epsilon},
\end{align*}
where 
\begin{equation*}
C_{c,c',\epsilon}=\underset{\underset{|v-u|\leq1}{|z| \leq A}}{\sup}|g(c(v-u)-c\epsilon z)-g(c'(v-u)-c'\epsilon z)|.
\end{equation*}
Let $Y_{c,c',\epsilon}$ be some centered gaussian variable with variance $C_{c,c',\epsilon}+\lambda^2 \ln(\frac{c'}{c})$. By using corollary \ref{cor:Kahane} of the appendix in a continuous version, we conclude that:
\begin{align*}
E[m_{c \epsilon}([0,c]^d)^p] & = E[(\int_{[0,c]^d}e^{X_{ c \epsilon}(x)-\frac{1}{2}E[X_{ c \epsilon}(x)^2]}dx)^p] \\
& = c^{dp} E[(\int_{[0,1]^d}e^{X_{ c \epsilon}(cu)-\frac{1}{2}E[X_{ c \epsilon}(cu)^2]}du)^p] \\
& \leq  c^{dp} E[(\int_{[0,1]^d}e^{X_{ c' \epsilon}(c'u)+Y_{c,c',\epsilon}-\frac{1}{2}E[(X_{ c' \epsilon}(c'u)+Y_{c,c',\epsilon})^2]}du)^p] \\
& = c^{dp}  (\frac{c'}{c})^{\frac{p(p-1)\lambda^2}{2}}e^{\frac{p(p-1)C_{c,c',\epsilon}}{2}}E[(\int_{[0,1]^d}e^{X_{ c' \epsilon}(c'u)-\frac{1}{2}E[X_{ c' \epsilon}(c'u)^2]}du)^p] \\
& = (\frac{c}{c'})^{dp-\frac{p(p-1)\lambda^2}{2}} e^{\frac{p(p-1)C_{c,c',\epsilon}}{2}} E[(\int_{[0,c']^d}e^{X_{ c' \epsilon}(x)-\frac{1}{2}E[X_{ c' \epsilon}(x)^2]}dx)  \\
& = (\frac{c}{c'})^{\zeta_p} e^{\frac{p(p-1)C_{c,c',\epsilon}}{2}} E[m_{c' \epsilon}([0,c']^d)^p] \\
\end{align*}
Taking the limit $\epsilon \rightarrow 0$ in the above inequality leads to:
\begin{equation}\label{eq:compsup}
 \frac{E[m([0,c]^d)^p]}{c^{\zeta_p}}  \leq  e^{\frac{p(p-1)C_{c,c'}}{2}} \frac{E[m([0,c']^d)^p]}{c'^{\zeta_p}}, 
 \end{equation}
where $C_{c,c'}=\underset{|v-u|\leq1}{\sup}|g(c(v-u))-g(c'(v-u))|$. Similarly, we have:
\begin{equation}\label{eq:compinf}
 \frac{E[m([0,c']^d)^p]}{c'^{\zeta_p}}  \leq  e^{\frac{p(p-1)C_{c,c'}}{2}} \frac{E[m([0,c]^d)^p]}{c^{\zeta_p}}.
 \end{equation}
Since $C_{c,c'}$ goes to $0$ as $c,c' \rightarrow 0$ ,we conclude by inequality (\ref{eq:compsup}) and (\ref{eq:compinf}) that $(\frac{E[m([0,c]^d)^p]}{c^{\zeta_p}})_{c >0}$ is a Cauchy sequence as $c \rightarrow 0$ bounded from below and above by positive constants. Therefore, there exists some $c_{p}>0$ such that:
\begin{equation*}
E[m([0,c]^d)^p] \underset{c \to 0}{\sim} c_{p} c^{\zeta_p}.
\end{equation*}
The same method can be applied to show that $\frac{c_{p}}{e^{\frac{p(p-1)g(0)}{2}}}$ is independent of $g$. This concludes the proof by setting $C_{p}=\frac{c_{p}}{e^{\frac{p(p-1)g(0)}{2}}}$.

\subsubsection{Proof of proposition \ref{prop:density}} 
We use the scaling relation (\ref{eq:invariance}) to compute the characteristic function of $m(B(0,c))$ for all $\xi$ in $\R$:
\begin{align*}
E[e^{i\xi m(B(0,c))}] & = E[e^{i\xi e^{\Omega_{c}}m(B(0,R))}] \\
& = E[\mathcal{F}(\xi m(B(0,R)))], \\
\end{align*}
where $\mathcal{F}$ is the characteristic function of $e^{\Omega_{c}}$. It is easy to show that for all $n \in \N$ there exists $C>0$ such that:
\begin{equation*}
|\mathcal{F}(\xi)| \leq \frac{C}{|\xi|^n}.
\end{equation*}
From this, we conclude by proposition \ref{prop:momneg} that:
\begin{equation*}
E[e^{i\xi m(B(0,c))}]  \leq  \frac{C}{|\xi|^n} E[\frac{1}{m(B(0,R))^n}] \leq  \frac{C'}{|\xi|^n}.
\end{equation*}
This implies the existence of a $C^{\infty}$ density.

\section{Appendix}

We give the following classical lemma first derived in \cite{cf:Kah}.

\begin{lemma}\label{lem:Kah} 
Let $(X_{i})_{1 \leq i \leq n}$ and $(Y_{i})_{1 \leq i \leq n}$ be two independent centered gaussian vectors and $(p_{i})_{1 \leq i \leq n}$ a sequence of positive numbers. If $\phi:\R_{+}\rightarrow \R$ is some smooth function with polynomial growth at infinity, we define:
\begin{equation*}
\varphi(t)=E[\phi(\sum_{i=1}^n p_{i}e^{Z_{i}(t)-\frac{1}{2}E[Z_{i}(t)^2]})],
\end{equation*}
with $Z_{i}(t)=\sqrt{t}X_{i}+\sqrt{1-t}Y_{i}$. Then, we have the following formula for the derivative:
\begin{equation}\label{eq:Kah}
\varphi'(t)=\frac{1}{2}\sum_{i=1}^n \sum_{j=1}^n p_{i}p_{j}(E[X_{i}X_{j}]-E[Y_{i}Y_{j}]) E[e^{Z_{i}(t)+Z_{j}(t)-\frac{1}{2}E[Z_{i}(t)^2]-\frac{1}{2}E[Z_{j}(t)^2]}\phi''(W_{n,t})],
\end{equation}
where:
\begin{equation*}
W_{n,t}=\sum_{k=1}^n p_{k}e^{Z_{k}(t)-\frac{1}{2}E[Z_{k}(t)^2]}
\end{equation*}
\end{lemma}

As a consequence of the above lemma, one can derive the following classical comparaison principle:

\begin{corollary}\label{cor:Kahane}
Let $(p_{i})_{1 \leq i \leq n}$ be a sequence of positive numbers. Consider $(X_{i})_{1 \leq i \leq n}$ and $(Y_{i})_{1 \leq i \leq n}$ two centered gaussian vectors such that:
\begin{equation*}
\forall i, j \quad E[X_{i}X_{j}] \leq E[Y_{i}Y_{j}].  
\end{equation*}
Then, for all convex function $F:\R\rightarrow \R_{+}$, we have:
\begin{equation}\label{eq:comparaison}
E[F(\sum_{i=1}^n p_{i}e^{X_{i}-\frac{1}{2}E[X_{i}^2]})], \leq E[F(\sum_{i=1}^n p_{i}e^{Y_{i}-\frac{1}{2}E[Y_{i}^2]})].
\end{equation}
\end{corollary}

We will also use another corollary:

\begin{corollary}\label{cor:Kah}
Let $(X_{i})_{1 \leq i \leq n}$ and $(Y_{i})_{1 \leq i \leq n}$ be two centered gaussian vectors such that:
\begin{itemize}
\item
$\forall i$, \quad \quad \quad $E[X_{i}^2]=E[Y_{i}^2]$.  
\item
$\forall i\not = j$, \quad $E[X_{i}X_{j}] \leq E[Y_{i}Y_{j}]$.  
\end{itemize}
Then, for all increasing function $F:\R\rightarrow \R_{+}$, we have:
\begin{equation}\label{eq:comp}
E[F(\sup_{1\leq i \leq n}Y_{i})] \leq E[F(\sup_{1\leq i \leq n}X_{i})].
\end{equation}
\end{corollary}

\proof
It is enough to show inequality (\ref{eq:comp}) for $F=1_{]x,+\infty[}$ for some $x \in \R$. Let $\beta$ be some positive parameter. Integrating equality (\ref{eq:Kah}) applied to the convex function $\phi: u \rightarrow 
e^{-e^{-\beta x}u}$ and the sequences $(\beta X_{i})$, $(\beta Y_{i})$, $p_{i}=e^{\frac{\beta^2}{2}E[X_{i}^2]}$, we get:
\begin{equation*}
E[e^{-\sum_{i=1}^{n}e^{\beta (X_{i}-x)}}] \leq E[e^{-\sum_{i=1}^{n}e^{\beta (Y_{i}-x)}}] 
\end{equation*}  
By letting $\beta \rightarrow \infty$, we conclude:
\begin{equation*}
P( \sup_{1\leq i \leq n}X_{i} < x) \leq P( \sup_{1\leq i \leq n}Y_{i} < x).
\end{equation*}  

\qed



\bigskip


\begin{thebibliography}{20}



\bibitem{cf:BaDeMu} Bacry, E., Delour, J., and Muzy, J.F.: Multifractal random walks, 
 \emph{Phys. Rev. E}, \textbf{64} (2001), 026103-026106.
\bibitem{cf:BaKoMu} Bacry E., Kozhemyak, A., Muzy J.-F.: Continuous
  cascade models for asset returns, available at
  www.cmap.polytechnique.fr/~bacry/biblio.html, to appear in \emph{Journal of Economic Dynamics and Control}. 
\bibitem{cf:BaMu} Bacry, E. and Muzy, J.F.: Log-infinitely divisible multifractal process, 
 \emph{Communications in Mathematical Physics}, \textbf{236} (2003), 449-475.


\bibitem{cf:Cha} Chainais, P.: Multidimensional infinitely divisible cascades. Application to the modelling of intermittency in turbulence, \emph{European Physical Journal B}, \textbf{51} no. 2 (2006), pp. 229-243.


\bibitem{cf:Cizeau} Cizeau, P., Gopikrishnan, P., Liu, Y., Meyer, M., Peng, C.K., Stanley, E.: Statistical properties of the volatility of price fluctuations, \emph{Physical Review E}, \textbf{60} no.2 (1999), 1390-1400.


\bibitem{cf:Co} Cont, R.: Empirical properties of asset returns: stylized facts and statistical issues, 
 \emph{Quantitative Finance}, \textbf{1} no.2 (2001), 223-236.

\bibitem{cf:DaVe} Daley D.J., Vere-Jones D.: An introduction to the
  theory of point processes, Springer-Verlag, (1988).
\bibitem{cf:DuRoVa} Duchon, J., Robert, R., Vargas, V.: Forecasting volatility with the multifractal random walk model, submitted to \emph{Mathematical Finance}, available at http://arxiv.org/abs/0801.4220.    


\bibitem{cf:Fr} Frisch, U.: \emph{Turbulence}, Cambridge University Press (1995).                                                                                                                                                                                                                                                                                                                                                                                                                                                                                                                                                                                                                                                                                                                                                              

\bibitem{cf:Gn} Gneiting, T.: Criteria of Polya type for radial positive definite functions, \emph{Proceedings of the American Mathematical Society}, \textbf{129} no. 8 (2001), 2309-2318.

\bibitem{cf:GrMacMa} Gray, Mathews, Macrobert: \emph{Bessel Functions}, Macmillan and co. (1922). 
\bibitem{cf:Kah} Kahane, J.-P.: Sur le chaos multiplicatif,
  \emph{Ann. Sci. Math. Qu{\'e}bec}, \textbf{9} no.2 (1985), 105-150.
\bibitem{cf:Kol} Kolmogorov A.N.: A refinement of previous hypotheses
  concerning the local structure of turbulence,
  \emph{J. Fluid. Mech.}, \textbf{13} (1962), 83-85. 

\bibitem{cf:LeTa} Ledoux, M., Talagrand, M.: \emph{Probability in Banach Spaces}, Springer-Verlag (1991). 
\bibitem{cf:Man} Mandelbrot B.B.: A possible refinement of the
  lognormal hypothesis concerning the distribution of energy in
  intermittent turbulence, \emph{Statistical Models and Turbulence},
  La Jolla, CA, Lecture Notes in Phys. no. 12, Springer, (1972), 333-351.
\bibitem{cf:Mol} Molchan, G. M.: Scaling exponents and multifractal dimensions for independent random cascades, \emph{Communications in Mathematical Physics}, \textbf{179} (1996), 681-702. 


\bibitem{cf:Obu} Obukhov A.M.: Some specific features of atmospheric
  turbulence, \emph{J. Fluid. Mech.}, \textbf{13} (1962), 77-81. 
\bibitem{cf:PaYu} Pasenchenko, O. Yu.: Sufficient conditions for the characteristic function of a two-dimensional isotropic distribution, \emph{Theory Probab. Math. Statist.}, \textbf{53} (1996), 149-152.
\bibitem{cf:Sch} Schwartz, L.: \emph{Théorie des distributions}, Hermann (1951).





\end{thebibliography}
\end{document}